\documentclass[12pt, a4paper,reqno]{amsart}
\usepackage[top=25mm, bottom=25mm, left=17mm, right=17mm]{geometry}
\usepackage{tikz}
\usepackage{nicematrix}
\usepackage{mathrsfs}
\usepackage{upgreek}

\newtheorem{theorem}{Theorem}[section]
\newtheorem{lemma}[theorem]{Lemma}
\newtheorem{proposition}[theorem]{Proposition}

\numberwithin{equation}{section}
\numberwithin{table}{section}
\numberwithin{figure}{section}
\newcommand{\abs}[1]{\lvert#1\rvert}

\newcommand{\wt}[1]{\widetilde{#1}}
\def\ontop#1#2{\setbox0\hbox{#2}\copy0\llap{\raise\ht0\hbox{#1}}}

\newcommand{\vertiii}[1]{{\left\vert\kern-0.25ex\left\vert\kern-0.25ex\left\vert #1 
	\right\vert\kern-0.25ex\right\vert\kern-0.25ex\right\vert}}
\newcommand{\eup}[1]{\langle#1\rangle}
\newcommand{\floor}[1]{\lfloor #1 \rfloor}
\newcommand{\ul}[1]{\underline#1}
\newcommand{\ol}[1]{\overline#1}
\newcommand{\bigabs}[1]{\Big\lvert#1\Big\rvert}

\begin{document}
\title{On the stability of the $L^{2}$ projection and the quasiinterpolant in the space of smooth
periodic splines.}
\author{D.C. Antonopoulos}
\address{Department of Mathematics, University of Athens, 15784 Zographou, Greece, and Institute of 
Applied and Computational Mathematics, FORTH, 70013 Heraklion, Greece}
\email{antonod@math.uoa.gr}
\author{V.A. Dougalis}
\address{Department of Mathematics, University of Athens, 15784 Zographou, Greece, and Institute of 
Applied and Computational Mathematics, FORTH, 70013 Heraklion, Greece}
\email{doug@math.uoa.gr}
\subjclass[AMS subject classification]{65D07,65M12}
\keywords{Smooth periodic polynomial splines, $L^{2}$ projection, Quasiinterpolant,
stability estimates, cyclic matrices, inverse of Gram matrix.}
\maketitle
\markright{\MakeUppercase{On  the Stability of the $L^{2}$ projection and the quasiinterpolant}}
\begin{abstract}
In this paper we derive stability estimates in $L^{2}$- and $L^{\infty}$- based Sobolev spaces for 
the $L^{2}$ projection and a family of quasiinterolants in the space of smooth, 1-periodic, 
polynomial splines defined on a uniform mesh in $[0,1]$. As a result of the assumed periodicity and 
the uniform mesh, cyclic matrix techniques and suitable decay estimates of the elements of the 
inverse of a Gram matrix associated with the standard basis of the space of splines, are used to 
establish the stability results. 
\end{abstract}
\maketitle
\section{Introduction}
In this paper we derive stability estimates in several norms for the $L^{2}$ projection and a family
of quasiinterpolants in the space of smooth, periodic, polynomial splines on a uniform mesh in a 
finite interval. Such estimates are of course useful in deriving approximation properties for the
$L^{2}$ projection and the quasiinterpolant. They have also been used in deriving optimal-order-of
-accuracy error estimates for Galerkin finite-element approximations to the solution of e.g. first-
order hyperbolic problems and nonlinear dispersive wave pde's, cf. e.g. \cite{bdk},
\cite{dk}, \cite{bdkmc}, \cite{adm}, \cite{ad1}, \cite{adk}. Stability estimates for the $L^{2}$ 
projection onto finite element spaces have been derived in several works. For example, see 
\cite{ddw} for estimates in $L^{\infty}$ in the case of a quasiuniform mesh in one dimension, 
\cite{ct} for estimates in $L^{p}$ and the Sobolev spaces $W_{\infty}^{p}$, $1\leq p\leq \infty$, in one and two 
dimensions for more general than quasiuniform types of meshes, and also the references of these 
papers. Here, we take advantage of the uniform mesh, and the periodicity of the underlying space of 
smooth splines, and use matrix methods to derive suitable decay estimates of the elements of the
inverse of a Gram matrix associated with the standard basis of the space of splines. To this end, we
use properties of cyclic (circulant) matrices listed in \cite{se}, and relevant results of Demko
\emph{et al.}, \cite{dms}, and Bini and Capovani, \cite{bc}. These decay properties enable us to
prove stability estimates for the $L^{2}$ projection in the $L^{2}$-based Sobolev spaces 
$H_{per}^{l}$ of periodic functions, for $l=0,\dots,r-1$, where $r$ is the order of the spline space,
and for continuous, periodic functions in $W_{l}^{\infty}$, for $l=0,\dots,r-1$. \par
Quasiinterpolants of continuous, periodic functions in the space of smooth, periodic splines have
been studied, among other, in \cite{s}, \cite{tw}, \cite{t}, \cite{ds}, and in  references of 
these works. All these quasiinterpolants achieve optimal-order accuracy in $L^{2}$ provided the
functions they approximate are smooth enough. Of particular interest is the Thom{\'e}e-Wendroff
quasiinterpolant, \cite{tw}, for which $L^{2}$ inner products of truncation errors with elements of
a special basis of the spline space are superaccurate due to cancellations. This enables one to prove
optimal-order $L^{2}$-error estimates for Galerkin approximations of the solutions of the periodic
initial-value-problems for the types of pde's previously mentioned.\par
In section 2 of the paper at hand we introduce the spline spaces and their standard basis. In 
section 3 we list a series of properties of cyclic matrices that will be needed in the sequel, 
mostly following \cite{se}. In section 4, based on results from \cite{dms} and \cite{bc}, we 
establish the required decay estimates for the elements of the inverse of the Gram matrix. These are
need in section 5 in order to establish stability results for the $L^{2}$ projection onto the spline
spaces in the Sobolev spaces $H_{per}^{l}$ and $W_{\infty}^{l}\cap H_{per}^{l}$ for 
$l=0,1,\dots r-1$. Finally, in section 6, we consider the quasiinterpolants and prove stability 
estimates for them in $H_{per}^{l}$ and $W_{\infty}^{l}\cap H_{per}^{l}$, for $l=1,\dots,r-1$, and in
$C_{per}$. These estimates do not depend of course on the decay results of section 4.\par
We use the following notation: For integer $k\geq 0$, $C_{per}^{k}$ ($C_{per}\equiv C_{per}^{0}$)
will denote the space of continuous, 1-periodic functions that are $k$ times continuously
differentiable. (By $C_{per}^{k}([0,1])$ we mean the restriction on $[0,1]$ of such functions).
Analogously, $L_{per}^{2}$ will denote the space of 1-periodic functions that are 
square-integrable over one period. The $L^{2}$ inner product, resp. norm, on $[0,1]$ will be denoted by
$(\cdot,\cdot)$, resp. $\|\cdot\|$. For integer $l\geq 0$, the norm on $[0,1]$ of the 
$L^{2}$-based Sobolev spaces of 1-periodic functions will be denoted by $\|\cdot\|_{l}$, while
the analogous norm of $W_{\infty}^{l}$ by $\|\cdot\|_{l,\infty}$. 
The Euclidean inner product, resp. norm, on $\mathbb{R}^{N}$  will be denoted by 
$\eup{\cdot,\cdot}$, resp. $\abs{\cdot}$, while
$\mathbb{P}_{r}$ will stand for the polynomials of degree at most $r$.
Finally, for $a>0$ we will  denote by $\floor{a}$ the largest integer that is less or equal to 
$a$.
\section{Spline spaces, basis} Let $r$ and $N$ be integers such that $r\geq 2$, $N\geq 4r$. Let
$h=1/N$ and $x_{i}=ih$, $i=0,1,\dots,N$, be a uniform partition of $[0,1]$. We consider the 
$N$-dimensional space of the 1-periodic, smooth, piecewise polynomial splines
\[
\mathcal{S}_{h}^{r}=\{v\in C_{per}^{r-2}[0,1] : v_{|(x_{i-1},x_{i})}\in
\mathbb{P}_{r-1}, \,\, i=1,2,\dots,N\},
\]
and the space
\[
\mathcal{S}_{h}^{1} = \{v\in L_{per}^{2}(0,1) : v_{|(x_{i-1},x_{i})}=\text{constant},
\,\, i=1,2,\dots,N\}.
\]
Following e.g. \cite{tw}, we define a standard basis of $\mathcal{S}_{h}^{r}$ as follows. Let 
$v_{1}$ be the characteristic function of $[-\tfrac{1}{2},\tfrac{1}{2}]$ and $v_{r}$ the convolution
of $v_{1}$ with itself $r-1$ times. Thus, $v_{r}$, with support $[-\tfrac{1}{2}r,\tfrac{1}{2}r]$, is
the $B$-spline of order $r$. If
\[
\phi_{l}(x)=v_{r}(h^{-1}x - l - \tfrac{r-2}{2}),\,\, l\in \mathbb{Z},\quad
\text{and}\quad
\Phi_{j}(x) =\sum_{l\in \mathbb{Z}}\phi_{j+lN}(x),\,\, j=1,2,\dots, N,
\]
then, the restrictions of $\{\Phi_{j}\}_{j=1}^{N}$ on $[0,1]$ are a basis of $\mathcal{S}_{h}^{r}$. \par
{\bf{Remarks:}} 1. Here we define $\Phi_{j}$ in terms of 
$\phi_{l}(x)=v_{r}(h^{-1}x - l - \tfrac{r-2}{2})$ and not in the standard way via 
$\phi_{l}(x)=v_{r}(h^{-1}x - l)$, since we wish that the support of $\Phi_{j}$ in $[0,1]$ to be
an interval or union of intervals with endpoints that are integer multiples of $h$. Thus, in 
$[0,1]$, supp$\Phi_{j}=[x_{j-1},x_{j+r-1}]$, for $j=1,2,\dots,N-(r-1)$, and 
supp$\Phi_{j}=[x_{j-1},1]\cup [0,x_{j-(N-(r-1))}]$ for $j=N-(r-2),N-(r-3),\dots,N$. \\
2. If $v\in\mathcal{S}_{h}^{r}$, then $v$ may be written in the form 
\[
v(x) = \sum_{j=1}^{N}V_{j}\Phi_{j}(x), \quad \text{or}\quad 
v(x)=\sum_{l\in\mathbb{Z}}V_{l}\phi_{l}(x),
\]
where the coefficients $V_{l}$ are periodic with period $N$, i.e. $V_{l}=V_{l+N}$, $l\in\mathbb{Z}$.
In addition it holds that $\Phi_{l}(x)=\Phi_{l+N}(x)$ for any integer $l$.\\
3. A basis of $\mathcal{S}_{h}^{1}$ is, obviously, the set $\{\Phi_{j}\}_{j=1}^{N}$, where
$\Phi_{j}(x)=v_{1}(h^{-1}x-j+\tfrac{1}{2})$.\\
4. If  $\{\Phi_{j}\}_{j=1}^{N}$, $\{\Psi_{j}\}_{j=1}^{N}$ are the bases of $\mathcal{S}_{h}^{r}$,
$\mathcal{S}_{h}^{r-1}$, respectively, then $\Phi_{j}'(x)=h^{-1}(\Psi_{j}(x) - \Psi_{j+1}(x))$.
Indeed, it follows by the definition of $v_{r}$ that
\[
v_{r}(x)=\int_{x-\frac{1}{2}}^{x+\frac{1}{2}}v_{r-1}(y)dy,
\]
from which $v_{r}'(x) = v_{r-1}(x+\frac{1}{2}) - v_{r-1}(x - \frac{1}{2})$, and since
$\Phi_{j}(x)=\sum_{l\in\mathbb{Z}}v_{r}(h^{-1}x - lN - j - \frac{r-2}{2})$, it follows that
\begin{align*}
\Phi_{j}'(x) & = h^{-1}\sum_{l\in\mathbb{Z}}
\Bigl[v_{r-1}\Bigl(h^{-1}x - lN-j-\frac{r-2}{2} +
\frac{1}{2}\Bigr) - v_{r-1}\Bigl(h^{-1}x -
lN-j-\frac{r-2}{2} -\frac{1}{2}\Bigr)\Bigr]\\
& = h^{-1}\sum_{l\in\mathbb{Z}}\Bigl[v_{r-1}\Bigl(h^{-1}x - lN-j - \frac{r-3}{2}\Bigr)
- v_{r-1}\Bigl(h^{-1}x - lN - j-1 - \frac{r-3}{2}\Bigr)\Bigr]\\
& = h^{-1}\bigl(\Psi_{j}(x) - \Psi_{j+1}(x)\bigr).
\end{align*}
5. If $G$  is the Gram matrix with elements $h^{-1}(\Phi_{k},\Phi_{j})$, then, cf. 
\cite[Lemma 2.1]{tw}, $G$ is a symmetric, positive definite, cyclic $N\times N$ matrix with
eigenvalues $g(2\pi l/N)$, $l=1,2,\dots,N$, where
\[
g(x) = \sum_{l\in\mathbb{Z}}\hat{v}_{r}(x + 2\pi l)^{2}, \qquad
\hat{v}_{r}(x)=\Biggl(\frac{2\sin\frac{1}{2}x}{x}\Biggr)^{r}.
\]
Moreover,
\[
\underline{g}\abs{V}^{2} \leq \eup{GV,V} \leq \overline{g}\abs{V}^{2},
\]
provided $\underline{g}\leq g(x)\leq \overline{g}$. It follows that: \\
(i)\,\, The elements of $G$ are numbers independent of $h$.\\
(ii)\, The function $g$ is periodic with period $2\pi$, $g(0)=g(2\pi)=1$, and
$g(\theta)\in [\underline{g},\overline{g}]$ for all $\theta\in\mathbb{R}$, 
\indent\,\,\,\, where $\overline{g}=1$
and $\underline{g} >0$. In addition, the maximum eigenvalue of $G$ is $\lambda_{\max}(G)=1$.\\
(iii)\, The eigenvalues of $G$ are included in the interval $[\underline{g},1]$ and 
$\underline{g}$ does not depend on $N$ but only \indent\,\,\,\,\, on $r$.\\
(iv)\,\, For each element of the inverse of $G$ there holds 
$\abs{(G^{-1})_{ij}}\leq 1/\lambda_{\min}(G)\leq 1/\underline{g}$, where $\lambda_{\min}(G)$ 
\indent \,\,\,\,\,
is the smallest eigenvalue of $G$. In particular, $\underline{g}\leq (G^{-1})_{11}\leq 1$.\\
(v)\,\, If $v_{h}=\sum_{j=1}^{N}V_{j}\Phi_{j}$ is an element of $\mathcal{S}_{h}^{r}$, then
\begin{equation}
\underline{g}h\abs{V}^{2} \leq \|v_{h}\|^{2} \leq \overline{g}h\abs{V}^{2}.
\label{eq21}
\end{equation} 
It is also well known that the following inverse inequality holds in $\mathcal{S}_{h}^{r}$: There
exists a constant $C$ independent of $h$, such that
\begin{equation}
\|v_{hx}\|\leq Ch^{-1}\|v_{h}\|,
\label{eq22}
\end{equation}
for all $v_{h}\in \mathcal{S}_{h}^{r}$.
\section{Cyclic matrices} 
In what follows we state some facts about cyclic (or circulant) matrices that will be useful in the
sequel. A $N\times N$ matric $C$ is called cyclic when
\begin{equation*}
c_{jk}=
\begin{cases}
c_{k-j+1}\,, &  \text{if} \,\,\, k\geq j,\\
c_{k-j+N+1}\,,& \text{if} \,\,\, k<j,
\end{cases}
\end{equation*}
where $c_{j}=c_{1j}$ are the elements of the first row of $C$. A cyclic matrix may be written in the
form 
\[
C = \sum_{j=1}^{N}c_{j}P^{j-1},
\]
where $P$ is the $N\times N$ permutation matrix whose columns are the vectors $e_{j}$ of the
standard basis of $\mathbb{R}^{N}$ written in the sequence
\begin{equation}
P=[e_{N},e_{1},e_{2},\dots,e_{N-1}].
\label{eq31}
\end{equation}
The following are properties of cyclic matrices, cf. e.g. \cite{se}:\\
(i)\,\,\, The matrix $P$ is cyclic.\\
(ii)\,\, $P^{N}=I_{N}$, where $I_{N}$ is the $N\times N$ identity matrix.\\
(iii)\, $PP^{T} = P^{T}P = I_{N}$. This follows from the fact that 
$PP^{T} = [Pe_{2},Pe_{3},\dots,Pe_{N},Pe_{1}]$, and \\ \indent\,\,\,\,\, then 
$PP^{T}=[e_{1},e_{2},\dots,e_{N-1},e_{N}]=I_{N}$. \\
(iv)\, The product of cyclic $N\times N$ matrices is cyclic.\\
(v)\,\,\, If $A$, $B$ are $N\times N$ cyclic matrices, then $AB=BA$.\\
(vi)\, The inverse of a cyclic matrix, if it exists, is cyclic.\\
(vii)\, If $C$ is a symmetric cyclic $N\times N$ matrix, then $c_{N+1-j}=c_{j+1}$, for 
$j=1,2,\dots,\floor{\frac{N-1}{2}}$. Such \indent \,\,\,\,\, a matrix may be written in the form
\begin{equation}
C = c_{1}I_{N} + \sum_{j=2}^{\frac{N}{2}}c_{j}\Bigl(P^{j-1} + P^{1-j}\Bigr) +
c_{\frac{N}{2} + 1}P^{\frac{N}{2}}, \quad \text{if}\,\,\, N\,\,\text{is even},
\label{eq32}
\end{equation}\indent\,\,\,\,
or
\begin{equation}
C = c_{1}I_{N} + \sum_{j=2}^{\frac{N+1}{2}}c_{j}\Bigl(P^{j-1} + P^{1-j}\Bigr),
\quad \text{if}\,\,\, N\,\,\text{is odd}.
\label{eq33}
\end{equation}\indent\,\,\,\,
To check this note that if $N$ is even, then 
\[
C = \sum_{j=1}^{N}c_{j}P^{j-1}=c_{1}I_{N} + \sum_{j=2}^{\frac{N}{2}}c_{j}P^{j-1}
+c_{\frac{N}{2}+1}P^{\frac{N}{2}} + \sum_{j=\frac{N}{2}+2}^{N}c_{j}P^{j-1}.
\]\indent\,\,\,\,
But $c_{j} = c_{N+2-j}$, $j=N/2+2, N/2+3,\dots,N$, from which we obtain 
\[
\sum_{j=\frac{N}{2}+2}^{N}c_{j}P^{j-1}=\sum_{j=\frac{N}{2}+2}^{N}c_{N+2-j}P^{j-1}
=\sum_{j=2}^{\frac{N}{2}}c_{j}P^{N-(j-1)}=\sum_{j=2}^{\frac{N}{2}}c_{j}P^{1-j},
\]\indent\,\,\,\,
and \eqref{eq32} follows. If $N$ is odd, then
\[
C = \sum_{j=1}^{N}c_{j}P^{j-1}=c_{1}I_{N} + \sum_{j=2}^{\frac{N+1}{2}}c_{j}P^{j-1}
+ \sum_{j=\frac{N+1}{2}+1}^{N}c_{j}P^{j-1}.
\]\indent\,\,\,\,
Here $c_{j}=c_{N+2-j}$, $j=(N+1)/2+1, (N+1)/2 +2,\dots,N$, giving 
\[
\sum_{j=\frac{N+1}{2}+1}^{N}c_{j}P^{j-1}= \sum_{j=\frac{N+1}{2}+1}^{N}c_{N+2-j}P^{j-1}
=  \sum_{j=2}^{\frac{N+1}{2}}c_{j}P^{N+2-j-1}= \sum_{j=2}^{\frac{N+1}{2}}c_{j}P^{1-j},
\]\indent\,\,\,\,
from which \eqref{eq33} follows.\\
(viii)\, If $C$ is an antisymmetric cyclic $N\times N$ matrix, then $c_{N+1-j}=-c_{j+1}$,
for $j=1,2,\dots,\floor{\frac{N-1}{2}}$, \indent\,\,\,\,\,\,\, with $c_{\frac{N}{2}+1}=0$, if $N$ 
is even. The matrix $C$ is written in the form 
\begin{equation}
C = \sum_{j=2}^{\floor{\frac{N+1}{2}}}c_{j}\Bigl(P^{j-1} - P^{1-j}\Bigr).
\label{eq34}
\end{equation}
(ix)\,\, If $V\in \mathbb{R}^{N}$ with $V=(V_{1},V_{2},\dots,V_{N})^{T}$ and $V_{l}=V_{l+N}$, 
$l\in\mathbb{Z}$, then $(P^{k}V)_{i}=V_{i+k}$, if \\ \indent\,\,\,\, $i=1,2\dots,N$, and $k$ 
any integer.
\section{An estimate of the elements of the inverse of $G$}
For the purposes of proving stability properties for the $L^{2}$ projection operator onto
$\mathcal{S}_{h}^{r}$, we shall prove some estimates for the elements of the inverse of the Gram
matrix $G$, introduced in Remark 5 of section 2.\par
Let $\{\Phi_{j}\}_{j=1}^{N}$ be the basis of $\mathcal{S}_{h}^{r}$ defined in section 2, let $G$ be
the matrix with elements $G_{ij}=h^{-1}(\Phi_{j},\Phi_{i})$, and $\Gamma$ be the inverse of $G$. The 
matrices $G$, $\Gamma$ are cyclic, symmetric, and positive definite. If 
$g_{j}=h^{-1}(\Phi_{j},\Phi_{1})$, then the first row of $G$ is the row vector
\[ 
(g_{1},g_{2},\dots,g_{r},0,\dots,0,g_{r},g_{r-1},\dots,g_{2}),
\]
and according to \eqref{eq32}, and
\eqref{eq33}, $G$ may be written as 
\begin{equation}
G = g_{1}I_{N} + \sum_{j=2}^{r}g_{j}\Bigl(P^{j-1} + P^{1-j}\Bigr).
\label{eq41}
\end{equation}
$\Gamma$ is a full matrix and its first row is
\begin{align*}
\gamma & =(\gamma_{1},\gamma_{2},\dots,\gamma_{\frac{N}{2}},\gamma_{\frac{N}{2}+1},
\gamma_{\frac{N}{2}}, \gamma_{\frac{N}{2}-1},\dots,\gamma_{2}), \quad \text{if}
\quad N \,\, \text{is even, or}\\
\gamma & =(\gamma_{1},\gamma_{2},\dots,\gamma_{\frac{N-1}{2}},\gamma_{\frac{N+1}{2}},
\gamma_{\frac{N+1}{2}}, \gamma_{\frac{N-1}{2}},\dots,\gamma_{2}), \quad \text{if}
\quad N \,\, \text{is odd}.
\end{align*}
Therefore, in view again of \eqref{eq32}, \eqref{eq33}, $\Gamma$ may be written in the form
\begin{equation}
\begin{aligned}
\Gamma & =\gamma_{1}I_{N} + \sum_{j=2}^{\frac{N}{2}}\gamma_{j}\Bigl(P^{j-1}+P^{1-j}\Bigr)
+ \gamma_{\frac{N}{2}+1}P^{\frac{N}{2}}, \quad \text{if}\quad N \,\, \text{is even},\\
\Gamma & =\gamma_{1}I_{N} + \sum_{j=2}^{\frac{N+1}{2}}\gamma_{j}\Bigl(P^{j-1}+P^{1-j}\Bigr),
\quad \text{if}\quad N \,\, \text{is odd}.
\end{aligned}
\label{eq42}
\end{equation}
\par
Demko \emph{et al.}, \cite[Proposition 2.2]{dms}, proved that if a banded matrix is positive
definite, then the elements of its inverse decay exponentially as they move away from the diagonal.
This result may be stated as follows for a symmetric, positive definite, banded matrix.
\begin{proposition} Let $B$ be a $N\times N$ positive definite, symmetric, banded matrix, for which
$B_{ij}=0$, if $\abs{i-j} > k$ for some positive integer $k$. If $\lambda_{\min}(B)$, 
$\lambda_{\max}(B)$ are the minimum and the maximum eigenvalue, respectively, of $B$ and 
$\lambda_{\mu}=(\lambda_{\min}(B))^{1/2}$, $\lambda_{m}=(\lambda_{\max}(B))^{1/2}$, then
\begin{equation}
\abs{(B^{-1})_{ij}} \leq C_{B} q_{B}^{-\abs{i-j}},
\label{eq43}
\end{equation}
where
\begin{equation}
C_{B} = \frac{1}{\lambda_{\mu}^{2}}\max \Bigl(1, \frac{(\lambda_{\mu} 
+ \lambda_{m})^{2}}{2\lambda_{m}^{2}}\Bigr),\quad
q_{B} = \Bigl(\frac{\lambda_{m}+\lambda_{\mu}}{\lambda_{m} - \lambda_{\mu}}\Bigr)^{1/k}.
\label{eq44}
\end{equation}
\end{proposition} 
Based on this result we will prove the following lemma for the inverse of $G$.
\begin{lemma} Let $\{\Phi_{j}\}_{j=1}^{N}$ be the basis of $\mathcal{S}_{h}^{r}$ defined in section 2,
$G$ be the $N\times N$ matrix defined by $G_{ij} = h^{-1}(\Phi_{j},\Phi_{i})$, and let $\Gamma$ be
the inverse of $G$. If $\gamma$ is the first row of $\Gamma$, there exist positive constants
$C_{1}$, $C_{2}$, $q$, independent of $N$, such that 
\begin{equation}
\abs{\gamma_{i}} \leq C_{1}q^{-(i-1)} + C_{2}q^{-(N-i)},
\label{eq45}
\end{equation}
for $i=r,r+1,\dots, \floor{N/2}+1$. Moreover, there exists a constant $C_{3}$, independent of $N$,
such that 
\begin{equation}
\sum_{i=r}^{\floor{\frac{N}{2}}+1}(1 + i)\abs{\gamma_{i}} \leq C_{3}.
\label{eq46}
\end{equation}
\end{lemma}
\begin{proof} Let $\wt{G}$ be the $N\times N$ symmetric, banded matrix with elements
$\wt{G}_{ij} = G_{ij}$, if $\abs{i-j}\leq r-1$, and $\wt{G}_{ij}=0$, if $\abs{i-j}>r-1$. Let
$\mathcal{G}$ be the $(N+2r-2)\times (N+2r-2)$ cyclic matrix with elements 
$\mathcal{G}_{ij}=\mathfrak{h}^{-1}(\Upphi_{j},\Upphi_{i})$, $1\leq i,j\leq N+2r-2$,
where $\mathfrak{h}=1/(N+2r-2)$, and $\{ \Upphi_{j} \}_{j=1}^{N+2r-2}$ is the basis of 
$\mathcal{S}_{\mathfrak{h}}^{r}$. Hence $\mathcal{G}$ is a
`cyclic extension' of $\wt{G}$, in the sense that $\wt{G}$ is obtained from $\mathcal{G}$ if we omit
the first $r-1$ and the last $r-1$ columns of $\mathcal{G}$, and also the first $r-1$ and the last
$r-1$ rows of $\mathcal{G}$. Following Bini and Capovani, cf. \cite[Proposition 4.2]{bc},
we obtain that $\lambda_{\min}(\mathcal{G})\leq \lambda_{\min}(\wt{G})$. Hence 
$\lambda_{\min}(\wt{G}) \geq \underline{g} >0$, i.e. $\wt{G}$ is positive definite, and by
Proposition 4.1, we have 
\begin{equation}
\abs{(\wt{G}^{-1})_{ij}} \leq C_{\wt{G}}q_{\wt{G}}^{-\abs{i-j}},
\label{eq47}
\end{equation}
where $C_{\wt{G}}$, $q_{\wt{G}}$ are defined as in Proposition 4.1 for $B=\wt{G}$ and $k=r-1$. Since
$\lambda_{min}(\wt{G})\geq\underline{g}>0$, it follows that $C_{\wt{G}}$ is bounded above by a
constant independent of $N$ (see Remark 5 (iii) in section 2.) In addition, using Proposition 4.2 of
\cite{bc}, we also conclude that $\lambda_{min}(\wt{G}) < \lambda_{max}(\wt{G})$ since the
eigenvalues of $\mathcal{G}$ have the same property. We conclude that $q_{\wt{G}}$ is also
independent of $N$. (Of course $C_{\wt{G}}$, $q_{\wt{G}}$ depend on $r$.) Now, the cyclic matrix
$G$ is written in the form $G=\wt{G}+\wt{W}$, where
\begin{equation*}
\wt{W} = 
\begin{pNiceMatrix}
0 &  & W \\
\hline
 & 0 & \\ 
\hline
W^{T} & & 0 
\end{pNiceMatrix}
, \qquad
W =
\begin{pmatrix}
g_{r} & g_{r-1} & \cdots & g_{2} \\
0 & g_{r} & \cdots & g_{3} \\
\vdots  & \ddots  & \ddots & \vdots  \\
 0 & \cdots & 0 & g_{r}
\end{pmatrix},
\end{equation*}
and therefore  $\wt{W}\gamma=(W\ul{\gamma},0,\dots,0,W^{T}\ol{\gamma})^{T}$, where
$\ul{\gamma} = (\gamma_{r},\gamma_{r-1},\dots,\gamma_{2})^{T}$,
$\ol{\gamma}=(\gamma_{1},\gamma_{2},\dots,\gamma_{r-1})^{T}$. Since
$G\gamma=e_{1}$ we have $\gamma = (\wt{G})^{-1}e_{1} - (\wt{G}^{-1})\wt{W}\gamma$, from which, for
$r\leq i\leq \floor{N/2}+1$, we see that
\begin{equation}
\begin{aligned}
\gamma_{i} & =(\wt{G}^{-1})_{i1} - \sum_{j=1}^{r-1}(W\ul{\gamma})_{j}(\wt{G}^{-1})_{ij}
-\sum_{j=1}^{r-1}(W^{T}\ol{\gamma})_{j}(\wt{G}^{-1})_{i,N-(r-1-j)} 
= (1 - (W\ul{\gamma})_{1})(\wt{G}^{-1})_{i1} \\
& \hspace{53pt} 
- (W^{T}\ol{\gamma})_{1}(\wt{G}^{-1})_{i,N-(r-2)} 
-\sum_{j=2}^{r-1}\Bigl((W\ul{\gamma})_{j}(\wt{G}^{-1})_{ij}
+ (W^{T}\ol{\gamma})_{j}(\wt{G}^{-1})_{i,N-(r-1-j)}\Bigr).
\end{aligned}
\label{eq48}
\end{equation}
But $g_{1}\gamma_{1} + 2(g_{2}\gamma_{2} + \dots g_{r}\gamma_{r})=1$, which yields
\[
(W\ul{\gamma})_{1} = \sum_{j=2}^{r}g_{j}\gamma_{j}=\frac{1}{2}(1-g_{1}\gamma_{1}).
\]
In addition, $(W^{T}\ol{\gamma})_{1}=g_{r}\gamma_{1}$, and therefore \eqref{eq48} may be written as
\begin{equation}
\gamma_{i}=\frac{1+g_{1}\gamma_{1}}{2}(\wt{G}^{-1})_{i1}
- g_{r}\gamma_{1}(\wt{G}^{-1})_{i,N-(r-2)}
-\sum_{j=2}^{r-1}\Bigl((W\ul{\gamma})_{j}(\wt{G}^{-1})_{ij}
+(W^{T}\ol{\gamma})_{j}(\wt{G}^{-1})_{i,N-(r-1-j)}\Bigr).
\label{eq49}
\end{equation}
Now, since $g_{1}+2(g_{2}+g_{3}+\dots+g_{r})=1$, $g_{j}>0$, for
$j=1,2,\dots,r$, $0<\gamma_{1}\leq 1$ and $\abs{\gamma_{j}}\leq 1/\ul{g}$,
it follows that 
$\abs{(W\ul{\gamma})_{j}} \leq 1/\ul{g}$ and $\abs{(W^{T}\ol{\gamma})_{j}} \leq 1/\ul{g}$,
for $j=2,3,\dots,r$. Hence, \eqref{eq48} and \eqref{eq43} give, for
$r\leq i\leq \floor{N/2}+1$,
\begin{align*}
\abs{\gamma_{i}} & \leq Cq^{-(i-1)}+  Cq^{-(N-r+2-i)}
+\frac{C}{\ul{g}}\Bigl(\sum_{j=2}^{r-1}q^{-(i-j)}
+ \sum_{j=2}^{r-1}q^{-(N-r+1+j-i)}\Bigr)\\
& = Cq^{-(i-1)} + Cq^{r-2}\cdot q^{-(N-i)} + \frac{C}{\ul{g}}\Bigl(
q^{-(i-1)}\cdot\frac{q^{r-1}-q}{q-1} + q^{-(N-i)}\cdot\frac{q^{r-2}-1}{q-1}\Bigr),
\end{align*}
where $C$ is a multiple of ${C}_{\wt{G}}$ by a constant and where, for simplicity, we have put 
$q = q_{\wt{G}}$. This finally yields \eqref{eq45} with
\[
C_{1} = \Bigl(1 + \frac{q^{r-1}-q}{\ul{g}(q-1)}\Bigr)C,\quad
C_{2} = \Bigl(q^{r-2} + \frac{q^{r-2}-1}{\ul{g}(q-1)}\Bigr)C.
\]
In order to prove \eqref{eq46} we consider the polynomial 
$p(x)=x^{r} + x^{r+1} + \dots +x^{M-1}=(x^{M}-x^{r})/(x-1)$, for $0<x\ne 1$, and
$M=\floor{N/2}+2$, and put $f(x)=xp(x)$. Then, it follows from \eqref{eq45} that
\[
\sum_{i=r}^{M-1}(1+i)\abs{\gamma_{i}}\leq C_{1}qf'(q^{-1})
+ C_{2}q^{-N}f'(q).
\]
Since 
\begin{align*}
f'(x)=\Bigl(\frac{x^{M+1}-x^{r+1}}{x-1}\Bigr)' & =
\frac{\bigl((M+1)x^{M}-(r+1)x^{r}\bigr)(x-1)-(x^{M+1}-x^{r+1})}{(x-1)^{2}}\\
& = \frac{Mx^{M+1}-rx^{r+1}-(M+1)x^{M}+(r+1)x^{r}}{(x-1)^{2}}\\
& < \frac{Mx^{M+1} + (r+1)x^{r}}{(x-1)^{2}},
\end{align*}
we get the estimate
\[
\sum_{i=r}^{M-1}(1+i)\abs{\gamma_{i}} \leq
C_{1}\frac{Mq^{-M} + (r+1)q^{-(r-1)}}{(q^{-1}-1)^{2}}
+ C_{2}\frac{Mq^{-N+M+1} + (r+1)q^{-N+r}}{(q-1)^{2}}.
\]
From this, the fact that 
\[
Mq^{-N+M+1} \leq \Bigl(\frac{N}{2}+2\Bigr)q^{-(\frac{N}{2}+2)} \cdot q^{5},
\]
and taking into account that the function $\omega(x)=xq^{-x}$ is bounded for $x>0$, we obtain
\eqref{eq46}.
\end{proof}
{\bf{Remarks:}} 1. It Follows from \cite[Theorem 1]{se}, that the elements of  $\Gamma$ for $r=2$
may be computed exactly.\\
2. From \eqref{eq46} and the fact that $\abs{\gamma_{i}}\leq 1/\ul{g}$, $1\leq i<r$, it follows that
there exists a constant $C_{4}$, independent of $N$, such that 
\begin{equation}
\sum_{i=1}^{\floor{\frac{N}{2}}+1}(1+i)\abs{\gamma_{i}}\leq C_{4}.
\label{eq410}
\end{equation} 
\section{Stability of the $L^{2}$ projection onto $\mathcal{S}_{h}^{r}$}
In this section we shall show the $L^{2}$-projection operator onto $\mathcal{S}_{h}^{r}$ is stable
in $H_{per}^{l}$ and $W_{\infty}^{l}$ for $l=0,1,\dots,r-1$.
\begin{proposition} Let $P_{h}$ be the $L^{2}$-projection operator onto $\mathcal{S}_{h}^{r}$. Then,
for $l=0,1,\dots,r-1$, there exist constants $C_{1}$, $C_{2}$, depending only on $l$ and $r$, such
that \\
(i)\,\, $\|\partial_{x}^{l}(P_{h}u)\|\leq C_{1} \|\partial_{x}^{l}u\|$, \quad
if \quad $u \in H^{l}_{per}$, and \\
(ii)\,\, $\|\partial_{x}^{l}(P_{h}u)\|_{\infty}\leq C_{2}\|\partial_{x}^{l}u\|_{\infty}$,
\quad if \quad $u\in W_{\infty}^{l}\cap H^{l}_{per}$.
\end{proposition}
\begin{proof} (i)\,\, If $l=0$ the estimate is obvious with $C_{1}=1$. Let now $1\leq l\leq r-1$ and
let $\{\Phi_{j}\}_{j=1}^{N}$ be the basis of $\mathcal{S}_{h}^{r}$, defined in section 2. Let
$P_{h}u=\sum_{j=1}^{N}c_{j}\Phi_{j}$. Then, applying Remark 4 of section 2, we get by induction that
\begin{equation}
\partial_{x}^{l}(P_{h}u) = h^{-l}\sum_{j=1}^{N}
\bigl((I_{N}-P^{-1})^{l}c\bigr)_{j}\Psi_{j},
\label{eq51}
\end{equation}
where $P$ is the `translation' matrix defined in \eqref{eq31}, $\{\Psi_{j}\}_{j=1}^{N}$ is the basis
of $\mathcal{S}_{h}^{r-l}$ and 
\[
c=(c_{1},c_{2},\dots,c_{N})^{T}
\]
is the solution of the linear system
$hGc=b$, with $G_{i,j}=h^{-1}(\Phi_{j},\Phi_{i})$, $1\leq i,j\leq N$ as usual, and
$b=(b_{1},b_{2},\dots,b_{N})^{T}$, $b_{i}=(u,\Phi_{i})$, $1\leq i\leq N$. We easily see by 
induction that 
\[
(I_{N} - P^{-1})^{l} = \sum_{m=0}^{l}\binom{l}{m}(-1)^{m}P^{-m}.
\]
Therefore, we get $hG\bigl((I_{N}-P^{-1})^{l}\bigr)c=(I_{N}-P^{-1})^{l}b=:\beta$, where
\begin{equation}
\beta_{i} = \sum_{m=0}^{l}\binom{l}{m}(-1)^{m}b_{i-m}.
\label{eq52}
\end{equation}
Now 
\begin{align*}
b_{i-m}=(u,\Phi_{i-m}) & = \int_{0}^{1}u(x)\Phi_{i}(x+mh)dx =
\int_{0}^{1}u(x-mh)\Phi_{i}(x)dx \\
& = \int_{0}^{1}u(x-mh)\Phi_{i-\wt{l}}(x-\wt{l}h)dx\\
& = \int_{0}^{1}u(x+(\wt{l}-m)h)\Phi_{i-\wt{l}}(x)dx
=:(\wt{u}_{l,m},\Phi_{i-\wt{l}}),
\end{align*}
where $\wt{l} = \floor{\frac{l}{2}}$ and $\wt{u}_{l,m}(x)=u(x+(\wt{l}-m)h)$. Hence
\begin{align*}
\wt{u}_{l,m}(x) & =\sum_{k=0}^{l-1}\frac{(\wt{l}-m)^{k}h^{k}}{k!}
(\partial_{x}^{k}u)(x) + \frac{1}{(l-1)!}\int_{\delta_{x}}(\partial_{x}^{l}u)(y)
\bigl(x+(\wt{l}-m)h - y\bigr)^{l-1}dy \\
& =: U_{l,m}(x) + \frac{1}{(l-1)!}J_{l,m}(x),
\end{align*}
where $\delta_{x}$ is the interval with endpoints $x$, $x+(\wt{l}-m)h$. It follows from 
\eqref{eq52} that
\begin{align*}
\beta_{i} & = \sum_{m=0}^{l}\binom{l}{m}(-1)^{m}(U_{l,m},\Phi_{i-\wt{l}}) +
\frac{1}{(l-1)!}\sum_{m=0}^{l}\binom{l}{m}(-1)^{m}(J_{l,m}, \Phi_{i-\wt{l}})\\
& = \sum_{k=0}^{l-1}\frac{h^{k}}{k!}(\partial_{x}^{k}u,\Phi_{i-\wt{l}})
\sum_{m=0}^{l}\binom{l}{m}(\wt{l}-m)^{k}(-1)^{m}
+ \frac{1}{(l-1)!}\sum_{m=0}^{l}\binom{l}{m}(-1)^{m}(J_{l,m},\Phi_{i-\wt{l}}).
\end{align*}
According to Lemma 5.2 that we prove at the end of the section, the sum 
$\sum_{m=0}^{l}\binom{l}{m}(\wt{l}-m)^{k}(-1)^{m}$, for $0\leq k\leq l$, is equal to zero. This
yields 
\begin{equation}
\beta_{i} = \frac{1}{(l-1)!}
\sum_{m=0}^{l}\binom{l}{m}(-1)^{m}(J_{l,m},\Phi_{i-\wt{l}}).
\label{eq53}
\end{equation}
Since
\begin{equation}
\abs{J_{l,m}(x)}^{2} \leq
\bigabs{\int_{\delta_{x}}(\partial_{x}^{l}u)^{2}(y)dy}\cdot
\bigabs{\int_{\delta_{x}}(x+(\wt{l}-m)h-y)^{2l-2}dy},
\label{eq54}
\end{equation}
and
\begin{align*}
\bigabs{\int_{\delta_{x}}(x+(\wt{l}-m)h-y)^{2l-2}dy}
& = \int_{0}^{\abs{\wt{l}-m}h}y^{2l-2}dy
=\frac{1}{2l-1}\abs{\wt{l}-m}^{2l-1}h^{2l-1} \\
& \leq \frac{1}{2l-1}\floor{\tfrac{l+1}{2}}^{2l-1}h^{2l-1}=:\wt{C}_{1,l}h^{2l-1},
\end{align*}
it follows that
\begin{align*}
\|J_{l,m}\|^{2} & \leq \wt{C}_{1,l}h^{2l-1}\sum_{j=1}^{N}
\int_{x_{j-1}}^{x_{j}}\bigabs{\int_{\delta_{x}}(\partial_{x}^{l}u)^{2}(y)dy}dx\\
& \leq \wt{C}_{1,l}h^{2l-1}\sum_{j=1}^{N}\int_{x_{j-1}}^{x_{j}}\int_{I_{l,m}}
(\partial_{x}^{l}u)^{2}(y)dydx,
\end{align*}
where $I_{l,m}$ is an interval of length $(\abs{\wt{l}-m}+1)h$. Therefore
\[
\|J_{l,m}\|^{2}\leq \wt{C}_{1,l}h^{2l}(\abs{\wt{l}-m}+1)\|\partial_{x}^{l}u\|^{2},
\]
for some constant $\wt{C}_{1,l}$. Finally
\begin{equation}
\|J_{l,m}\| \leq C_{1,l}h^{l}\|\partial_{x}^{l}u\|, \quad \text{if}\quad
u\in H^{l}_{per},
\label{eq55}
\end{equation}
for some constant $C_{1,l}$. Moreover
\[
\abs{J_{l,m}(x)} = \bigabs{\int_{\delta_{x}}(\partial_{x}^{l}u)(y)
\bigl(x+(\wt{l}-m)h - y\bigr)^{l-1}dy},
\]
and therefore 
\[
\abs{J_{l,m}(x)}\leq \|\partial_{x}^{l}u\|_{\infty}\int_{0}^{\abs{\wt{l}-m}h}y^{l-1}dy,
\]
giving
\begin{equation}
\|J_{l,m}\|_{\infty} \leq C_{2,l} h^{l}\|\partial_{x}^{l}u\|_{\infty}, \quad
\text{if}\quad u\in H^{l}_{per}\cap W_{l}^{\infty},
\label{eq56}
\end{equation}
for some constant $C_{2,l}$. From \eqref{eq53}, taking into account \eqref{eq55} and \eqref{eq21} we
see that 
\begin{align*}
\abs{\beta}^{2} & = \frac{1}{(l-1)!}
\sum_{m=0}^{l}\binom{l}{m}(-1)^{m}
\bigl(J_{l,m},\sum_{i=1}^{N}\beta_{i}\Phi_{i-\wt{l}}\bigr)\\
& \leq C_{l}\|J_{l,m}\|h^{0.5}\abs{\beta}\leq C_{3,l}h^{l+0.5}
\|\partial_{x}^{l}u\|\abs{\beta},
\end{align*}
where $C_{l}$ , $C_{3,l}$ are constants depending on $l$. Hence
\[
\abs{\beta} \leq C_{3,l}h^{l+0.5}\|\partial_{x}^{l}u\|.
\]
From this estimate, the fact that $hG(I_{N}-P^{-1})^{l}c=\beta$, \eqref{eq21} and \eqref{eq51}, it
follows that there exists a constant $C=C(l,r)$, such that
\[
\|\partial_{x}^{l}(P_{h}u)\|\leq h^{-l+0.5}
\abs{(I_{N}-P^{-1})^{l}c} \leq Ch^{-l+0.5} \cdot h^{l-0.5}
\|\partial_{x}^{l}u\|,
\]
from which (i) follows with $C_{1}=C$. \par
(ii)\,\, Let $1\leq l\leq r-1$. If now $u\in H^{l}_{per}\cap W_{l}^{\infty}$, then from \eqref{eq53}
and \eqref{eq56}, we obtain
\[
\max_{1\leq i\leq N}\abs{\beta_{i}} \leq C_{4,l}h^{l+1}\|\partial_{x}^{l}u\|_{\infty},
\]
for some constant $C_{4,l}$. This estimate, the definition of $\beta$, and \eqref{eq410}, imply that
there is a constant $C=C(l,r)$, such that
\[
\max_{1\leq i\leq N}\abs{((I_{N}-P^{-1})c)_{i}}
\leq Ch^{l}\|\partial_{x}^{l}u\|_{\infty}.
\]
Therefore, from this inequality, \eqref{eq51}, and the fact that $\sum_{j=1}^{N}\Psi_{j}=1$, we see
that (ii) follows. \par
If $l=0$ and $\Gamma$ is the (cyclic, symmetric) inverse of $G$, then \eqref{eq410} implies that
$\|\Gamma\|_{\infty}\leq C_{4}$, and therefore 
$\abs{c}_{\infty}=\abs{\Gamma b}_{\infty}\leq C_{4}\abs{b}$. But $b_{i}=(u,\Phi_{i})$, so that
\[
\abs{b_{i}}\leq \|u\|_{\infty}\int_{I_{i}}\abs{\Phi_{i}(x)}dx, \quad 1\leq i\leq N,
\]
where $I_{i}$ is an interval (or union of intervals) of length $rh$. Therefore, there exists a
constant $C$, independent of $h$, such that 
\[
\abs{b}_{\infty}:=\max_{1\leq i\leq N}\abs{b_{i}}\leq Ch\|u\|_{\infty},
\]
i.e. that $\abs{c}_{\infty}\leq C_{4}\cdot Ch\|u\|_{\infty}$. Given that $0\leq \Phi_{j}(x)\leq 1$
for $j=1,2,\dots,N$ and all $x\in [0,1]$, it follows that
\[
\abs{(P_{h}u)(x)} \leq \sum_{j=1}^{N}\abs{c_{j}} \Phi_{j}(x)\le \sum_{j=1}^{N}\abs{c_{j}}
\leq C\|u\|_{\infty},
\]
i.e. $\|P_{h}u\|_{\infty}\leq C\|u\|_{\infty}$. The proof of proposition 5.1. is now complete.
\end{proof}
\begin{lemma} Let $l\geq 0$ be an integer, and for a nonnegative integer $k$ let
\[
a_{k} = \sum_{m=0}^{l}\binom{l}{m}m^{k}(-1)^{m}.
\]
Then, if $l=0$ and $k\geq 0$, or $l\geq 1$ and $0\leq k\leq l-1$, we have that $a_{k}=0$.
\end{lemma}
\begin{proof} For $l=0$, $k\geq 0$ we immediately get the conclusion. If $l\geq 1$, 
$0\leq k\leq l-1$, consider the polynomial 
\[
p(x) = (1+x)^{l} = \sum_{m=0}^{l}\binom{l}{m}x^{m}.
\]
Then
\[
p^{(k)}(x) = \frac{l!}{(l-k)!}(1+x)^{l-k} =
\sum_{m=k}^{l}\binom{l}{m}\frac{m!}{(m-k)!}x^{m-k}:=b_{k}(x),
\]
and therefore $p^{(k)}(-1) = 0 = b_{k}(-1)$, for $k=0,1,\dots,l-1$, and also $a_{0}=b_{0}(-1)=0$,
$a_{1}= - b_{1}(-1) = 0$. Assuming $a_{j}=0$ for $j=0,1,\dots,k-1$, we will show that $a_{k}=0$. We
have
\begin{align*}
b_{k}(-1) & = (-1)^{k}\sum_{m=k}^{l}\binom{l}{m}\frac{m!}{(m-k)!}(-1)^{m}\\
& =(-1)^{k}
\sum_{m=k}^{l}\binom{l}{m}m(m-1)\cdot\cdot\cdot\bigl(m-(k-1)\bigr)(-1)^{m}\\
& = (-1)^{k}\sum_{m=1}^{l}\binom{l}{m}\bigl(m^{k}+ q_{k-1}(m)\bigr)(-1)^{m},
\end{align*}
where $q_{k-1}(m)$ is a polynomial of $m$ of degree $k-1$. Since $b_{k}(-1)=0$ it follows that
\begin{equation}
\sum_{m=1}^{l}\binom{l}{m}\bigl(m^{k} + q_{k-1}(m)\bigr)(-1)^{m}=0.
\label{eq57}
\end{equation}
From our hypothesis that $a_{j}=0$ for $j=0,1,2,\dots,k-1$, we obtain
\[
\sum_{m=1}^{l}\binom{l}{m}q_{k-1}(m)(-1)^{m}=0,
\]
and, hence, using \eqref{eq57}, we get $a_{k}=0$, and the inductive step is complete.
\end{proof}
\section{Stability of a family of quasiinterpolants in $\mathcal{S}_{h}^{r}$}
We shall prove now that a family of quasiinterpolants in $\mathcal{S}_{h}^{r}$, cf. e.g.
\cite{s}, \cite{tw}, \cite{t}, \cite{ds}, is stable in $H_{per}^{l}$ for $1\leq l\leq r-1$, in
$(C_{per}, \|\cdot\|_{\infty})$, and in $W_{l}^{\infty}$ for $1\leq l\leq r-1$. We recall from the
above references that if $u$ is a continuous, 1-periodic function, then a \emph{quasiinterpolant}
in $\mathcal{S}_{h}^{r}$, is defined as 
\[
Q_{h}u= \sum_{j=1}^{N}u(jh)\wt{\Phi}_{j}\,,
\]
where $\{\wt{\Phi}_{j}\}_{j=1}^{N}$ is a modified basis of $\mathcal{S}_{h}^{r}$ defined as
$\wt{\Phi}_{j}=q_{0}\Phi_{j}+\sum_{k=1}^{r-1}q_{k}(\Phi_{j+k} + \Phi_{j-k})$, where
$\{\Phi_{j}\}_{j=1}^{N}$ is the standard basis of $\mathcal{S}_{h}^{r}$ as introduced in section 2.
Therefore, if $v\in \mathcal{S}_{h}^{r}$, then,
\[
v = \sum_{j=1}^{N}V_{j}\Phi_{j}=\sum_{j=1}^{N}\wt{V}_{j}\wt{\Phi}_{j},
\] 
with $V=Q\wt{V}$, where $Q$ is a symmetric $N\times N$ cyclic matrix with first row
\[
(q_{0},q_{1},\dots,q_{r-1},0,\dots,0,q_{r-1},\dots,q_{1}).
\]
Actually, see Remark 1 in section 2, the quasiinterpolant that we define 
here satisfies $(Q_h u)(x) = (Q^{TW}_{h} u)(x - h(r - 2)/2)$, but due to 
periodicity the results of this section are valid for $Q^{TW}_{h}$  as well.
 Note that for the Thom{\'e}e-Wendroff quasiinterpolant, \cite{tw}, the $q_{j}$'s are chosen so that the trigonometric polynomial 
\[
q(\xi)=\sum_{j=-(r-1)}^{r-1}q_{j}e^{-ij\xi}=
\sum_{j=0}^{r-1}\delta_{j}\Bigl(\sin\frac{\xi}{2}\Bigr)^{2j},
\]
corresponds to coefficients obtained from the series 
\[
\sum_{j=0}^{\infty}\delta_{j}\tau^{2j}=\Bigl(\frac{\arcsin\tau}{\tau}\Bigr)^{r}.
\]
As a result of this choice, this quasiinterpolant has several superconvergent approximation
properties which will not concern us here.
\begin{proposition} Let $Q_{h}$ be a quasiinterpolant in $\mathcal{S}_{h}^{r}$ defined as above.
Then: \\
(i)\,\, There exists a constant $C_{0}$ independent of $h$, such that 
$\|Q_{h}u\|_{\infty}\leq C_{0}\|u\|_{\infty}$, for all $u\in C_{per}$.\\
Moreover, for $l=1,2,\dots,r-1$, there exist constants $C_{1}$, $C_{2}$ depending only on $l$ and
$r$ such that \\
(ii)\,\, $\|\partial_{x}^{l}(Q_{h}u)\|_{\infty} \leq C_{1} \|\partial_{x}^{l}u\|_{\infty}$, \,\, if 
\,\, $u\in W_{\infty}^{l}\cap H^{l}_{per}$, and \\
(iii)\, $\|\partial_{x}^{l}(Q_{h}u)\| \leq C_{2}\|\partial_{x}^{l}u\|$,
\,\,if \,\, $u\in H_{per}^{l}$.
\end{proposition}
\begin{proof} (i)\,\, We have
\begin{equation}
Q_{h}u=\sum_{j=1}^{N}u(jh)\wt{\Phi}_{j} = \sum_{j=1}^{N}(Q\wt{u})_{j}\Phi_{j},
\label{eq61}
\end{equation}
where $Q = q_{0}I_{N} + \sum_{m=1}^{r-1}q_{m}(P^{m} + P^{-m})$ and 
$\wt{u}=(u_{1},u_{2},\dots,u_{N})^{T}$, $u_{j}=u(jh)$, $1\leq j\leq N$. Since
\begin{equation}
(Q\wt{u})_{j}=q_{0}u_{j} + \sum_{m=1}^{r-1}q_{m}(u_{j+m} + u_{j-m}),
\label{eq62}
\end{equation}
it follows that
\[
\max_{1\leq j\leq N}\abs{(Q\wt{u})_{j}} \leq 
(\abs{q_{0}}+2\sum_{m=1}^{r-1}\abs{q_{m}})\|u\|_{\infty} =: C_{0}\|u\|_{\infty},
\]
i.e. the result (i) of the proposition. \par
From \eqref{eq61} and Remark 4 of section 2,
\begin{equation}
\partial_{x}^{l}(Q_{h}u) =h^{-l}\sum_{j=1}^{N}\omega_{j}^{l}\Psi_{j},
\label{eq63}
\end{equation}
where $\{\Psi\}_{j=1}^{N}$ the basis of $\mathcal{S}_{h}^{r-l}$ and   
$\omega_{j}^{l} = \bigl[(I_{N} - P^{-1})^{l}(Q\wt{u})\bigr]_{j}$. Hence, from the relation
\[
(I_{N}-P^{-1})^{l}=\sum_{k=0}^{l}\binom{l}{k}(-1)^{k}P^{-k},
\]
and \eqref{eq62}, we obtain
\begin{equation}
\omega_{j}^{l} = q_{0}\sum_{k=0}^{l}\binom{l}{k}(-1)^{k}u_{j-k}
+ \sum_{m=1}^{r-1}q_{m}\sum_{k=0}^{l}\binom{l}{k}(-1)^{k}(u_{j-k+m} + u_{j-k-m}).
\label{eq64}
\end{equation}
Now, for integer $k$ it holds that 
\[
u_{j-k} = \sum_{\gamma=0}^{l-1}\frac{(-1)^{\gamma}k^{\gamma}h^{\gamma}}{\gamma !}
(\partial_{x}^{\gamma}u)(jh) + J(l,j,k),
\]
where
\[
J(l,j,k) = \frac{-1}{(l-1)!}\int_{(j-k)h}^{jh}
(\partial_{x}^{l}u)(x)\bigl((j-k)h - x\bigr)^{l-1}dx.
\] 
Therefore, from \eqref{eq64} and Lemma 5.2, we get that
\begin{equation}
\omega_{j}^{l}=q_{0}\sum_{k=0}^{l}{l}\binom{l}{k}(-1)^{k}J(l,j,k)
+ \sum_{m=1}^{r-1}q_{m}\sum_{k=0}^{l}\binom{l}{k}(-1)^{k}
\bigl(J(l,j,k-m) + J(l,j,k+m)\bigr).
\label{eq65}
\end{equation}
Since for integer $k\in \{-2(r-2),\dots,r-1\}$, it holds that
\[
\abs{J(l,j,k)} \leq \frac{\abs{k}^{l}h^{l}}{l!}\|\partial_{x}^{l}u\|_{\infty},
\]
and 
\[
\abs{J(l,j,k)}^{2}  
\leq \frac{\abs{k}^{2l-1}h^{2l-1}}{(2l-1)((l-1)!)^{2}}
\bigabs{\int_{(j-k)h}^{jh}(\partial_{x}^{l}u)^{2}(x)dx}
\leq \frac{\abs{k}^{2l-1}h^{2l-1}}{(2l-1)((l-1)!)^{2}}
\int_{(j-2r+4)h}^{(j+r-1)h}(\partial_{x}^{l}u)^{2}(x)dx,
\]
we obtain from \eqref{eq65} that
$\abs{\omega_{j}^{l}} \leq \wt{C}_{1}h^{l}\|\partial_{x}^{l}u\|_{\infty}$,
and
\[
\abs{\omega_{j}^{l}}^{2}\leq \wt{C}_{2}h^{2l-1}
\int_{(j-2r+4)h}^{(j+r-1)h}(\partial_{x}^{l}u)^{2}(x)dx.
\]
From these two inequalities, \eqref{eq63}, and \eqref{eq21}, we get the estimates (ii) and (iii)
of the Proposition.
\end{proof}
\bibliographystyle{amsalpha} 

\end{document}